\documentclass[12pt]{amsart}
\usepackage{amsxtra}

\theoremstyle{plain}

\newcommand{\cleqn}{\setcounter{equation}{0}}
\newcommand{\clth}{\setcounter{theorem}{0}}
\newcommand {\sectionnew}[1]{\section{#1}\cleqn\clth}

\newtheorem{theorem}{Theorem}[section]
\newtheorem{lemma}[theorem]{Lemma}
\newtheorem{definition-theorem}[theorem]{Definition-Theorem}
\newtheorem{proposition}[theorem]{Proposition}
\newtheorem{corollary}[theorem]{Corollary}
\newtheorem{definition}[theorem]{Definition}
\newtheorem{example}[theorem]{Example}
\newtheorem{remark}[theorem]{Remark}
\newtheorem{notation}[theorem]{Notation}
\newtheorem{assumption}[theorem]{Assumption}
\newtheorem{lemma-definition}[theorem]{Lemma-Definition}
\newcommand \bth[1] { \begin{theorem}\label{t#1} }
\newcommand \ble[1] { \begin{lemma}\label{l#1} }

\newcommand \bpr[1] { \begin{proposition}\label{p#1} }
\newcommand \bco[1] { \begin{corollary}\label{c#1} }
\newcommand \bde[1] { \begin{definition}\label{d#1}\rm }
\newcommand \bex[1] { \begin{example}\label{e#1}\rm }
\newcommand \bre[1] { \begin{remark}\label{r#1}\rm }

\newcommand \bnota[1] {\begin{notation}\label{n#1}\rm }
\newcommand \bas[1] { \begin{assumption}\label{a#1}\rm }

\newcommand {\eth} { \end{theorem} }
\newcommand {\ele} { \end{lemma} }

\newcommand {\epr} { \end{proposition} }
\newcommand {\eco} { \end{corollary} }
\newcommand {\ede} { \end{definition} }
\newcommand {\eex} { \end{example} }
\newcommand {\ere} { \end{remark} }
\newcommand {\enota} { \end{notation} }
\newcommand {\eas} {\end{assumption}}

\newcommand \thref[1]{Theorem \ref{t#1}}
\newcommand \leref[1]{Lemma \ref{l#1}}
\newcommand \prref[1]{Proposition \ref{p#1}}
\newcommand \coref[1]{Corollary \ref{c#1}}

\newcommand \exref[1]{Example \ref{e#1}}
\newcommand \reref[1]{Remark \ref{r#1}}
\newcommand \lb[1]{\label{#1}}

\def \Cset {{\mathbb C}}

\def \al {\alpha}

\def \la {\lambda}

\def \ra  {\rightarrow}           

\def \la {\langle}
\def \ra {\rangle}
            
\def \wt {\widetilde}

\def \Ad { {\mathrm{Ad}} }

\def \rank { {\mathrm{rank}} }

\def \g  {\mathfrak{g}}   
\def \h  {\mathfrak{h}}

\def \n  {\mathfrak{n}}

\def \b  {\mathfrak{b}}

\def \la{{\langle}}
\def \ra{{\rangle}}

\def \pist {\pi_{\rm st}}

\DeclareMathOperator \Aut { {\mathrm{Aut}} }

\newcommand \im { {\mathrm{im}} }

\renewcommand{\qed}{\begin{flushright} {\bf Q.E.D.}\ \ \ \ \
                  \end{flushright} }    

\def \hs {\hspace{.2in}}

\def \lara {\la \, \, , \, \ra}

\def \gog {\g \oplus \g}

\def \lara {\la  \, , \, \ra}

\def \hs {\hspace{.2in}}
\def \dw {\bar{w}}

\def \Th {\theta}
\def \th {\theta}
\def \pTh {\pi_\Th}

\def \piGs {\pi_{\scriptscriptstyle G^*}}

\def \gog {\g \oplus \g}
\def \Gdia {G_{\rm diag}}
\def \gdia {\g_{\rm diag}}
\def \gst {\g_{\rm st}^*}
\def \hwt {\h_{w\Th}}
\def \Twt {T_{w\Th}}
\def \Gso {G^*\cdot_\Th (h_0\dw)}
\def \cTh {\cdot_\Th}
\def \sc {{\scriptscriptstyle C}}
\def \Wc {W_\sc}
\def \Wcm {W_{\sc}^-}
\def \mc {m_{\sc}}
\def \rk {{\rm rank}}
\def \dk {\dim \ker}
\def \Aut {{\rm Aut}}

\def \th {\theta}
\def \pth {\pi_\th}

\def \tG {G^n}
\def \tth  {\tilde{\th}}
\def \ptth {\pi_{\tth}}
\def \sC {{\scriptscriptstyle C}}
\def \Gst {G^*}
\def \lrw {\longrightarrow}
\def \pist {\pi_{\rm st}}
\def \Pist {\Pi_{\rm st}}
\def \GuvC {G^{u,v}_{\sC}}

\def \SwC {S^{\dw}_{{\sC}}(h)}
\def \SwCz {S^{\dw}_{{\sC}}(h_0)}
\begin{document}
\setlength{\baselineskip}{1.2\baselineskip}
\title[On a Poisson structure on twisted conjugacy classes]
{On the $T$-leaves and the ranks of a Poisson structure on twisted conjugacy classes}
\author[Jiang-Hua Lu]{Jiang-Hua Lu}
\address{Jiang-Hua Lu,
Department of Mathematics,
Hong Kong University,
Pokfulam Rd., Hong Kong}
\email{jhlu@maths.hku.hk}
\date{}
\begin{abstract} Let $G$ be a connected complex semisimple Lie group with a fixed maximal torus $T$ and a Borel
subgroup $B \supset T$.  For an arbitrary automorphism $\theta$ of $G$, 
we introduce a holomorphic Poisson structure $\pi_\theta$ on $G$ which is invariant under the $\th$-twisted conjugation
by $T$ and has the property that 
every $\theta$-twisted conjugacy class of $G$ is a Poisson 
subvariety with respect to $\pi_\theta$. We describe the $T$-orbits of symplectic leaves, called $T$-leaves, of $\pi_\theta$ and
compute the dimensions of the symplectic leaves (i.e, the ranks) of $\pth$. We 
give the lowest rank of $\pth$ in any given $\theta$-twisted conjugacy class, and we
relate the lowest possible rank locus of $\pth$ in $G$ with spherical
$\th$-twisted conjugacy classes of $G$. In particular, we
show that 
$\pi_\theta$ vanishes somewhere on $G$ if and only if $\theta$ induces an involution on the Dynkin diagram of $G$, and that in such a
case a $\Th$-twisted conjugacy class $C$ contains a vanishing point of $\pi_\theta$ if and only if $C$ is spherical.
\end{abstract}
\maketitle

\sectionnew{Introduction  and statement of results}\lb{intro}

\subsection{Introduction}\lb{subsec-intro}
Let $G$ be a connected complex semi-simple Lie group and 
fix a Borel subgroup $B$ of $G$ and a maximal torus $T$ contained in $B$.
It is well-known that the choice of the pair $(B, T)$ 
gives rise to a so-called {\it standard} multiplicative holomorphic Poisson structure $\pist$ on $G$ that vanishes on $T$, and that the Poisson Lie group
$(G, \pist)$
is the semi-classical limit of the quantum group $\Cset_q[G]$, the Hopf algebra dual of the
quantized enveloping algebra $U_q \g$ of $G$ (see \cite{C-P:guide, Eting-Schiff}). The Poisson group $(G, \pist)$
is reviewed in $\S$\ref{subsec-piG}.

Let $\Th \in \Aut(G)$. Define the $\theta$-twisted conjugation of $G$ on itself by
\begin{equation}\lb{eq-th-1}
g_1 \cdot_\theta g = g_1g \theta(g_1)^{-1}, \hs g, g_1\in G,
\end{equation}
and call its orbits the $\Th$-twisted conjugacy classes in $G$.
Using a general construction from \cite{LiB-M} and 
\cite[Theorem 2.3]{L-Y:DQ}, we define
a holomorphic Poisson structure $\pTh$ on $G$ with the properties that every $\theta$-twisted 
conjugacy class $C$ in $G$ is a Poisson submanifold with respect to $\pTh$ and that $(C, \pTh)$ is  a 
Poisson homogeneous space \cite{dr:homog} of the Poisson Lie group $(G, \pist)$
with respect to the $\theta$-twisted
conjugation. 
The precise definition of $\pTh$ is given in $\S$\ref{subsec-dfn}. 

The Poisson structure 
$\pTh$ is invariant under the $\theta$-twisted conjugation by elements in $T$. If $\Sigma \subset G$ is a symplectic leaf of $\pTh$ in $G$, 
the set $T\cTh \Sigma = \cup_{h \in T} h\cTh \Sigma$ is called the $T$-orbit
of $\Sigma$ in $G$, or a $T$-leaf of $\pth$. 
When $\theta={\rm Id}_G$, the identity automorphism of $G$, the Poisson structure $\pTh$, in particular
its $T$-leaves and some $T$-equivariant Poisson resolutions of
certain Poisson subvarieties of $(G, \pTh)$, are studied in \cite{E-L:grothen}.

In this paper, for an arbitrary $\th \in {\rm Aut}(G)$, we describe the $T$-leaves of $\pth$ and
compute the dimensions of the symplectic 
leaves, also called the ranks, of $\pth$.  We describe the lowest rank of $\pth$ in any given 
$\theta$-twisted conjugacy class, and we relate the
lowest possible rank locus of  $\pth$ with spherical $\th$-twisted conjugacy classes in $G$. 
For the case of $\theta={\rm Id}_G$, 
our result on the $T$-leaves
of $\pth$ strengthens that given in \cite{E-L:grothen}. Work in this paper makes use of the results in \cite{CLT, Lu:formula}
on intersections of Bruhat cells and $\th$-twisted conjugacy classes and 
on the element
$\mc$ in the Weyl group of $G$ naturally associated to a $\th$-twisted conjugacy class $C$ in $G$. 
This paper can be regarded as first applications of results in \cite{Chan, CLT, Lu:formula}
to Poisson geometry (see \cite[Remark 3.10]{Lu:formula}). It can also be regarded as a sequel to \cite{E-L:grothen}.

\subsection{Statement of results}\lb{subsec-statements}
Recall that we fix, from the very beginning, a maximal torus $T$ of $G$ and a Borel subgroup $B$ containing $T$.
To simplify the statements on the Poisson structure $\pth$, 
we assume in this section that $\th(B) = B$ and $\th(T) = T$. 
The case of an arbitrary $\th \in {\rm Aut}(G)$ is treated in $\S$\ref{sec-th-arbitrary}.

Let $W=N_G(T)/T$ be the Weyl group, 
where $N_G(T)$ is the stabilizer of $T$ in $G$. Let $l: W \to {\mathbb N}$ be the length function on $W$
and let $\leq$ be the Bruhat order on $W$.
Let $B_-$ be the Borel subgroup of $G$ such that $B \cap B_- = T$. 
Then one has the Bruhat decompositions 
\[
G = \bigsqcup_{w \in W} B w B = \bigsqcup_{w \in W} BwB_-.
\]
where $BwB_-$ has co-dimension 
$l(w)$ in $G$ for $w \in W$. 
For a $\th$-twisted conjugacy class $C$ in $G$, let 
$\Wcm = \{w \in W: C \cap (BwB_-) \neq \emptyset\}.$
The
sets $\Wcm$ were first studied in \cite{CLT}. It particular, it is shown in \cite[$\S$2.4]{CLT} that  
\[
\Wcm = \{w \in W: w \leq \mc\},
\]
where $\mc$ is the unique element in $W$ such that $C \cap (B\mc B)$ is dense 
in $C$. 
The elements $\mc$ play an important role in the study of spherical conjugacy classes and their applications to the theory of
representations of quantum groups (see works of 
N. Cantarini, G. Carnovale, and M. Costantini in \cite{CCC, Car-2008, Car-2012, Co}). 
For $G$ simple and $\theta ={\rm Id}_G$, the list of all $\mc$'s, as $C$ runs over all conjugacy  classes in $G$, is 
given in \cite[$\S$3]{CLT} using results from \cite{CCC}. For $G$ of classical type and for every conjugacy class
$C$ in $G$, the element $\mc$ is explicitly computed in \cite{Chan}. For $G$ simple and $\theta$ an outer
automorphism, the list of the $\mc$'s, as $C$ runs over all $\theta$-twisted conjugacy  classes in $G$, is given in
\cite{Lu:formula}.

Our main results on the $T$-leaves of $\pth$ are summarized in the following 
\thref{th-T-orb-mc}, where for $w \in W$, $1 +w\th$ is regarded as a linear operator on the Lie algebra
$\h$ of $T$, with $1$ denoting the identity operator on $\h$.

\bth{th-T-orb-mc} The intersections $C \cap (BwB_-)$, where $C$ is a $\th$-twisted conjugacy class in $G$
and $w \in W$ is such that $w \leq \mc$, are precisely all the
$T$-leaves of $\pth$ in $G$, and the symplectic leaves in $C \cap (BwB_-)$ 
have dimension  equal to 
\[
\dim C - l(w) -\dim \ker (1+w\theta).
\]
\eth

A slightly stronger version of \thref{th-T-orb-mc}, where $T$ can be  replaced by certain subtori of $T$, is proved in \thref{th-TT-orb}.
Symplectic leaves in $C \cap (B w B_-)$ are described in \coref{co-leaves} as connected components of certain submanifolds 
of $C \cap (B w B_-)$.

Recall that a $\th$-twisted conjugacy class in $G$ is said to be {\it spherical} if it admits an open orbit under
the $\th$-twisted conjugation by $B$. Our second main result relates the lowest possible rank of $\pth$ in $G$ and the collection 
${\mathcal S}_\th$ of all 
spherical $\th$-twisted conjugacy classes in $G$. More precisely, for $g \in G$, let $\rk(\pth(g))$ be the rank of
$\pth$ at $g$, i.e., the 
dimension of the symplectic leaf of $\pth$ passing through $g$, and let $\rk(1-\th^2)$ be the rank of the linear operator $1-\th^2$ on
$\h$.

\bth{th-min-rank}
One has $\rk(\pth(g)) \geq \rk(1-\th^2)$ for every $g \in G$, and  
\[
\{g \in G: \rk(\pth(g)) = \rk(1-\th^2)\} = \bigsqcup_{C \in {\mathcal S}_\th} C \cap (B \mc B_-).
\]
In particular, the set $\{g \in G: \rk(\pth(g)) = \rk(1-\th^2)\}$ is nonempty if and only if spherical
$\th$-twisted conjugacy classes
in $G$ exist. 
\eth

Denote by $\Gamma$ both the Dynkin diagram of $G$ and the set of simple roots determined by
$(T, B)$. Since the induced action of $\th$ on $\h^*$ permutes the simple roots, it induces an automorphism on $\Gamma$
which is also denoted by $\th$. A consequence of \thref{th-min-rank} is the following conclusion on the
zero locus $Z(\pth) = \{g \in G: \pth(g) = 0\}$ of $\pth$.

\bco{co-0-locus}
1) $Z(\pth) \neq \emptyset$  if and only $\th^2 = 1 \in {\rm Aut}(\Gamma)$. 

2) When $\th^2 = 1 \in {\rm Aut}(\Gamma)$, one has
\[
Z(\pth)= \bigsqcup_{C \in {\mathcal S}_\th} C \cap (B \mc B_-),
\]
and $C \cap (B \mc B_-)$ is a single $T$-orbit for every $C \in {\mathcal S}_\th$.
\eco

Examples related to $SL(n+1, {\mathbb C})$ and the triality automorphism for $D_4$ are given in $\S$\ref{sec-T-leaves}
and $\S$\ref{sec-spherical}. In section $\S$\ref{sec-ptth} we study in some detail 
the 
Poisson structure $\ptth$ on $G^n$, where $n \geq 2$, and $\tth\in {\rm Aut}(G^n)$ is given by 
\[
\tth(g_1, g_2, \ldots, g_n) = (g_2, \ldots, g_n, g_1), \hs g_j \in G, 1 \leq j \leq n.
\]
In particular, for $n = 2$, we show that $(G \times G, \ptth)$ is isomorphic to the Drinfeld double $(G \times G, \Pist)$
(see $\S$\ref{subsec-piG}) 
of the Poisson Lie group $(G, \pist)$. Consequently,  for the diagonal action of $T$ on $G \times G$ by left (or right)
translations, the $T$-leaves of $\Pist$ in $G \times G$ are precisely the 
{\it double Bruhat cells
associated to conjugacy classes}, i.e., the submanifolds 
\[
G^{u, v}_\sC = \{(k_1, k_2) \in BuB \times B_- v B_-: \, k_1k_2^{-1} \in C\} \subset G \times G,
\]
where $u, v \in W$ and $C$ is any conjugacy class $C$ in $G$ (see also \reref{re-proof-pim} for a more direct proof of this fact).
Note that when $C = \{e\}$, $G^{u, v}_{\sC}$ is isomorphic to the double Bruhat cell $G^{u, v}= BuB \cap B_-v B_-$,
which
have been studied intensively \cite{b-f-z:III, f-z:double} in connection  with total positivity and cluster algebras.  
Further studies of double Bruhat cells
associated to conjugacy classes will be carried out elsewhere.
  
\subsection{Acknowledgment} The author would like to thank 
G. Carnovale, K. Y. Chan and V. Mouqiun for helpful discussions. She also thanks the referees for helpful comments
and particularly for \reref{re-referee}. This work is partially supported
by the Research Grants Council of the Hong Kong SAR, China (GRF HKU 703712P). 

\sectionnew{Definition of the Poisson structure $\pTh$ on $G$}\lb{sec-pth-def}

\subsection{The Poisson Lie group $(G, \pist)$ and its dual $(G^*, \piGs)$}\lb{subsec-piG} As in 
$\S$\ref{subsec-intro}, let $G$ be a connected complex semi-simple Lie group with Lie algebra $\g$ and fix
a pair $(B, T)$, where $B$ is a Borel subgroup of $G$ and $T \subset B$ a maximal torus of $G$.
Let $B_-$ be the Borel subgroup of $G$ such that $B \cap B_- = T$, and 
let $N$ and $N_-$ be respectively the uniradicals of
$B$ and $B_-$. The Lie algebras of $T, B, B_-, N$ and $N_-$ will be denoted by $\h, \b, \b_-, \n$ and $\n_-$ respectively.

Let $\lara_\g$ be  the Killing form of $\g$, and equip
 the direct product Lie algebra $\g \oplus \g$ with the symmetric
bilinear form 
\begin{equation}\lb{eq-lara-gog}
\la (x_1, y_1), \, (x_2, y_2) \ra_{\gog} = \la x_1, x_2 \ra_\g - \la y_1, y_2 \ra_\g, \hs x_1, x_2, y_1, y_2 \in \g.
\end{equation}
Let $\g_{\rm diag} = \{(x, x): x \in \g\}$ and let
\[
\g_{\rm st}^* = \{(x_+ + y, \, -y + x_-): \; x_+ \in \n, x_- \in \n_-, y \in \h\} \subset \gog.
\]
Then $\gog = \g_{\rm diag} + \g_{\rm st}^*$ is a direct sum of vector spaces and both $\g_{\rm diag}$ and $\g_{\rm st}^*$
are {\it Lagrangian subalgebras} of the quadratic Lie algebra $(\gog, \lara_{\gog})$ in the sense that they are Lie subalgebras
of $\gog$ and are 
maximal isotropic with
respect to $\lara_{\gog}$.  The decomposition $\g \oplus \g = \g_{\rm diag} + \g^*_{\rm st}$ is
called the 
{\it standard Lagrangian splitting} of $\gog$. 
The connected complex subgroups of $G \times G$ with Lie algebras $\gdia$
and $\g_{\rm st}^*$  are respectively $\Gdia = \{(g, g): g \in G\}$ and 
\begin{equation}\lb{eq-Gs}
G^* = \{(nh, mh^{-1}): n \in N, m \in N_-, h \in T\}.
\end{equation}
By the theory of Poisson Lie groups (\cite{C-P:guide, Eting-Schiff}, 
\cite[Appendix]{E-L:grothen}, \cite[Proposition 2.2]{L-Y:DQ}), the Lagrangian splitting $\gog = \gdia + \g_{\rm st}^*$
induces multiplicative holomorphic Poisson structures $\pist$ on $G$ and $\piGs$ on $G^*$, making
$(G, \pist)$ and $(G^*, \piGs)$ a dual pair of complex Poisson Lie groups. 
The $r$-matrix associated to the splitting $\gog = \gdia + \g_{\rm st}^*$ is the element $R \in \wedge^2(\gog)$ given by 
\begin{equation}\lb{Lambda}
R = \frac{1}{2} \sum_{j=1}^n (\xi_j \wedge x_j) \in \wedge^2 (\g \oplus \g), 
\end{equation}
where $\{x_j\}_{j = 1}^{n}$ is any basis  of $\gdia$ and 
$\{\xi_j\}_{j=1}^{n}$ the basis of $\gst$ such that $\la x_j, \xi_k\ra_{\gog} = \delta_{jk}$ for 
$1 \leq j, k \leq n$. 

The product group $G \times G$ carries the multiplicative Poisson structure 
\begin{equation}\lb{eq-pi-pm}
\Pist = R^r - R^l, 
\end{equation}
where 
$R^l$ and $R^r$ denote respectively the left and right invariant
bivector fields on $G \times G$ with value $R$ at the identity element. 
The Poisson Lie group $(G \times G, \Pist)$ is called a  {\it Drinfeld double} \cite{k-s:quantum}
of the Poisson Lie group $(G, \pist)$. 
It follows from the definition that $\Pist$ is invariant under the following action of $T \times T$ on $G \times G$:
\[
(h_1, h_2)\cdot (g_1,g_2) = (h_1g_1h_2^{-1},\, h_1 g_2h_2^{-1}), \hs h_1, h_2 \in T, \, g_1, g_2 \in G. 
\]
Moreover, the embeddings $(G, \pist)\cong (G_{\rm diag}, \pist) \hookrightarrow (G \times G, \Pist)$
and $(G^*, -\piGs) \hookrightarrow (G \times G, \Pist)$ are Poisson \cite[Proposition 2.2]{L-Y:DQ}.
We will return to the Poisson structure $\Pist$ on $G \times G$ in $\S$\ref{sec-ptth}.

\subsection{The definition of $\pTh$}\lb{subsec-dfn}
Let $\theta \in \Aut(G)$, and let $G_\theta = \{(g, \theta^{-1}(g)): g \in G\}$. 
Identify $(G \times G)/G_\theta$ with $G$ by the isomorphism
\[
\eta: \;\; (G \times G)/G_\theta \longrightarrow G: \; \; \; 
(g_1, g_2)G_\theta \longmapsto g_1 \theta(g_2)^{-1}, \hs g_1, g_2 \in G.
\]
The natural left action of $G \times G$ on $(G \times G)/G_\theta$ becomes that of $G \times G$ on $G$ by
\begin{equation}\lb{eq-GG-G}
(g_1, g_2)\cdot_\theta g = g_1 g \theta(g_2)^{-1}, \;\;\; g_1, g_2, g \in G.
\end{equation}
Note the action in  \eqref{eq-GG-G} restricted to $G_{\rm diag} \subset G \times G$ is the $\theta$-twisted conjugation of $G$
on itself given in \eqref{eq-th-1}, and, by abuse of notation, we are denoting both actions by $\cdot_\th$. 
Let $\kappa$ be the Lie algebra anti-homomorphism from $\gog$ to the Lie algebra of holomorphic
vector fields on $G$ induced by the action in \eqref{eq-GG-G}, i.e., 
\begin{equation}\lb{eq-kappa}
\kappa(x, y) = x^r - \theta(y)^l, \hs x, y \in \g, \;g \in G,
\end{equation}
where for $x \in \g$, $x^r$ and $x^l$ respectively denote the right and left invariant vector field on $G$ with value
$x$ at the identity element of $G$.  Define the bivector 
field $\pTh$ on $G$ by
\[
\pTh = (\kappa \wedge \kappa)(R),
\]
where $R \in \wedge^2 (\gog)$ is given in \eqref{Lambda}. 
The following \prref{pr-prop-pi-Th} is a special case of \cite[Proposition 2.2, Theorem 2.3]{L-Y:DQ}.

\bpr{pr-prop-pi-Th}
1) $\pTh$ is a Poisson bi-vector field on $G$, and every $\theta$-twisted conjugacy class
$C$ in $G$ is a Poisson submanifold of $(G, \pi_\th)$;

2) The twisted conjugation action
\[
(G, \pist) \times (G, \pTh) \longrightarrow (G, \pTh),\;\;\;
(g_1, g) \longmapsto g_1 \cdot_\th g = g_1 g \theta(g_{1})^{-1}, \hs g_1, g \in G,
\]
is a Poisson action of the Poisson Lie group $(G, \pist)$ on the Poisson manifold $(G, \pth)$;  

3) The symplectic leaves of $\pTh$ in $G$ are the connected components of the intersections of $\theta$-twisted conjugacy classes 
and $G^*$-orbits in $G$, where $G^*$ acts on $G$ 
as a subgroup of $G \times G$ via \eqref{eq-GG-G}.
\epr

\bre{re-DeConcini}
The $\theta$-twisted conjugacy class through the identity element of $G$ is isomorphic to $G/G^\th$, 
where $G^\th = \{g \in G: \th(g) = g\}$.
The Poisson manifold $(G/G^\theta, \pi_\th)$ is an example of a {\it De Concini} Poisson homogeneous space
of $(G, \pist)$ considered in \cite[Section 6.4]{En-Et}. The standard Lagrangian splitting $\gog = \g_{\rm diag} + \g_{\rm st}^*$
induces a Poisson structure $\Pi_0$ on the variety ${\mathcal L}$ of Lagrangian subalgebras of $(\gog, \lara_{\gog})$, which is considered in
\cite{E-L:cplx}.
When $G$ is of adjoint type, $G \cong (G \times G)/G_\th$ can be identified with 
the $(G \times G)$-orbit in ${\mathcal L}$ through the point $\g_\th = \{(x, \theta^{-1}(x)): x \in \g\} \in {\mathcal L}$,
and the induced Poisson structure on $G$ coincides with $\pTh$.
\ere

\bre{re-Th-Th} For $g \in G$, let $\Ad_g: G \to G, \Ad_g(g_1) = gg_1g^{-1}$ for $g_1 \in G$.
If $\theta, \theta^\prime \in {\rm Aut}(G)$ are such that 
$\theta= \Ad_{g} \theta^\prime$ for some $g \in G$, it is easy to see from the
definitions that the right translation on $G$ by $g$ intertwines the $\th$-twisted conjugation 
and the $\th^\prime$-twisted conjugation and is 
a Poisson isomorphism from $(G, \pth)$ to $(G, \pi_{\th^\prime})$. 
For an arbitrary $\theta \in {\rm Aut}(G)$, since all Borel subgroups of $G$ are conjugate to each other, there exists $g_1 \in G$
such that $\theta(B) = g_1Bg_1^{-1}$. Since $T$ and $g_1^{-1}\theta(T)g_1$ are two maximal tori of $G$ both contained in $B$, there
exists $b \in B$ such that $b^{-1}g_1^{-1}\theta(T)g_1b=T$. Let $g = g_1b$, and let $\theta^\prime = \Ad_{g^{-1}} \theta$. Then 
$\theta^\prime$ stabilizes both $B$ and $T$.  To study the geometry of the Poisson structure $\pth$ and in particular the 
$T$-leaves of $\pth$, we can thus
assume that $\theta(B) = B$ and $\th(T) = T$. This remark will be used in $\S$\ref{sec-th-arbitrary}.
\ere

Assuming that $\theta(B) = B$ and $\theta(T) = T$, we now give an explicit formula for $\pi_\theta$.

Let $\{y_i\}_{i = 1}^{k}$ be a basis of $\h$ such that
$2\la y_i, y_j\ra_\g = \delta_{ij}$ for $1 \leq i, j \leq k = \dim \h$. Let $\Delta_+ \subset \h^*$ be the set of positive roots
determined by $B$, and for each $\al \in \Delta_+$, let $E_\al$ and $E_{-\al}$ be root vectors for $\al$ and $-\al$ respectively such that
$\la E_\al, E_{-\al}\ra_\g = 1$. The following \leref{le-pth-explicit} in the case when $\th = {\rm Id}_G$ is given in
\cite[$\S$2.2]{E-L:grothen}.

\ble{le-pth-explicit}
One has
\begin{equation}\lb{eq-pth-explicit}
\pth=\sum_{i=1}^{k} \theta(y_i)^l \wedge y_i^r - \sum_{\alpha \in\Delta_+} E_{\alpha}^{r} \wedge \theta(E_{-\alpha})^l
+ \frac{1}{2} \sum_{\alpha \in \Delta_+}
\left(E_\alpha^r \wedge E_{-\alpha}^r + E_\alpha^l \wedge E_{-\alpha}^l\right).
\end{equation}
\ele

\begin{proof}
Using the bases  
\begin{align}\lb{eq-x-i}
\{x_i\}& = \{(y_1, y_1),  \;(y_2, y_2), \ldots, (y_k, y_k),  \;
(E_\alpha, \, E_\alpha), \;(E_{-\alpha}, E_{-\alpha}): 
\alpha \in \Delta_+\}\\
\lb{eq-xi-i}
\{\xi_i\} & = \{(y_1, -y_1),  \;(y_2, -y_2), \ldots, (y_k, -y_k),\;  
(0, -E_{-\alpha)},\; (E_{\alpha}, 0): 
\alpha \in \Delta_+\},
\end{align}
for $\gdia$ and $\gst$ and by \eqref{Lambda}, one has
\begin{align}\lb{eq-R}
R &= \sum_{i=1}^k (y_i, 0) \wedge (0, y_i) +\sum_{\al \in \Delta_+} (E_\al, 0) \wedge (0, E_{-\al}) \\
\nonumber 
& \;\; + \frac{1}{2}\sum_{\al \in \Delta_+} \left((E_\al, 0) \wedge (E_{-\al}, 0) +(0, E_{\al}) \wedge (0, E_{-\al})\right).
\end{align}
It follows from the definition of $\pi_\theta$ that 
\[
\pi_\theta = \sum_{i=1}^{k} \theta(y_i)^l \wedge y_i^r - \sum_{\alpha \in\Delta_+} E_{\alpha}^{r} \wedge \theta(E_{-\alpha})^l
+ \frac{1}{2} \sum_{\alpha \in \Delta_+}
\left(E_\alpha^r \wedge E_{-\alpha}^r + \theta(E_\alpha)^l \wedge \theta(E_{-\alpha})^l\right).
\]
Using the facts that $\theta$ preserves $\lara_\g$ and that 
the induced action of $\theta$ on $\h^*$ permutes the positive roots, one proves \eqref{eq-pth-explicit}.
\end{proof}
  
\sectionnew{$T$-leaves of $\pth$}\lb{sec-T-leaves}
\subsection{The intersections $C \cap (B wB_-)$}\lb{subsec-CBB}
Throughout $\S$\ref{sec-T-leaves}-$\S$\ref{sec-spherical}, we assume that 
$\theta(B)=B$ and $\theta(T)=T$. Then $\theta(B_-) = B_-$.  Recall the subgroup $G^* \subset B\times B_-$ 
given in \eqref{eq-Gs} and let $G^*$ act on $G$ via \eqref{eq-GG-G}. 
Since $\th(B) = B$ and $\th(B_-) = B_-$, every $G^*$-orbit in $G$ lies in $B wB_-$ for a unique $w \in W$. 
By \prref{pr-prop-pi-Th}, every symplectic leaf of $\pth$ in $G$ lies in the intersection $C \cap (B w B_-)$ for
a unique $\th$-twisted conjugacy class $C$ and a unique $w \in W$. One thus needs to first study the intersections
$C \cap (B w B_-)$ and in particular to know when such an intersection is nonempty. We now recall some results
from \cite{CLT}.

Let $C$ be an arbitrary $\th$-twisted conjugacy class in $G$.
Since $C$ is irreducible, there is a unique element
$\mc \in W$ such that $C \cap (B\mc B)$ is dense in $C$.  Let
\[
\Wc = \{w \in W: C \cap (BwB) \neq \emptyset\} \hs \mbox{and} \hs \Wcm= \{w \in W: C \cap (BwB_-) \neq \emptyset\}.
\]
It is clear that $\Wc \subset \Wcm$.  By \cite[$\S$1]{CCC},
$\mc$ is the unique maximal element of $\Wc$ with respect to the Bruhat order $\leq $ on $W$.
Since $\theta(T) = T$, one has $\theta(N_G(T)) = N_G(T)$ and thus an induced map 
\begin{equation}\lb{eq-th-W}
\th: \;\;\; W \longrightarrow W, \;\;\; \th(gT) = \th(g)T, \hs g \in N_G(T).
\end{equation}
Note that $\th(w)(h) = \th w(\th^{-1}(h))$ for $h \in T$. 
For $w \in W$, define the
$\theta$-twisted conjugacy class of $w$ in $W$ to be $\{vw\theta(v^{-1}): v \in W\}$.
It is shown in \cite[$\S$2.4]{CLT} that for every $\th$-twisted conjugacy class $C$ in $G$,
$\mc$ is the unique maximal length element in its $\theta$-twisted conjugacy class in $W$. Moreover \cite[Corollary 3.3 and Theorem 3.8]{CLT}, 
$m_{\scriptscriptstyle C}^2 = 1$ and  $\th(\mc) = \mc$. Thus
\begin{equation}\lb{eq-mc-inv}
(\mc \th)^2 = \th^2 \in {\rm End}(\h).
\end{equation}

\bpr{pr-mc} \cite[$\S$2.4]{CLT} One has $\Wcm = \{w \in W: \; w \leq \mc\}$ for every $\theta$-twisted conjugacy class $C$ in $G$.
\epr

\bex{ex-SL}
Let $G = SL(n+1, {\mathbb C})$ and $\th = {\rm Id}_G$. Let $B$ and $B_-$ be the subgroups of $G$ consisting, respectively,
of all upper-triangular
and lower triangular matrices.
Identify the Weyl group of $G$ with the symmetric group 
$S_{n+1}$ on $n+1$ letters. For a conjugacy class $C$ in $SL(n+1, {\mathbb C})$, let 
$r(C)= \min\{\rank(g-cI): c \in {\mathbb C}\}$
for any  $g \in C$, and let $l(C) = \min (r(C), [(n+1)/2])$, where 
$[(n+1)/2]=m$ if $n = 2m$ or $n = 2m-1$.  Let $m_0 = 1$, and for an integer
$1 \leq l \leq [(n+1)/2]$, let $m_l \in S_{n+1}$ be the involution with the cycle decomposition
\[
m_l = (1, n+1) (2, n) \cdots (l,\; n+2-l).
\]
By \cite[Corollary 4.3]{CLT}, $\mc = m_{l({\scriptscriptstyle C})}$
for any conjugacy class $C$ in $SL(n+1, {\mathbb C})$.
This example will be continued in \exref{ex-SL-2}.
\eex

\ble{le-irre}
For any $\th$-twisted conjugacy class $C$ in $G$ and any $w \in W$ such that $w \leq \mc$, the
intersection $C \cap (BwB_-)$ is a smooth and connected submanifold of
$G$ and has dimension equal to $\dim C - l(w)$.
\ele

\begin{proof}  
Since $\Gdia \cap (B \times B_-)$ is connected, it follows from
\cite[Corollary 1.5]{Ri} that the
intersection $C \cap (BwB_-)$ is transversal and irreducible, so $C \cap (BwB_-)$ is smooth and connected,  and 
\[
\dim (C \cap (BwB_-)) = \dim C + \dim (BwB_-) - \dim G = \dim C - l(w).
\]
\end{proof}

\subsection{Symplectic leaves in $C \cap (B w B_-)$}\lb{subsec-leaves-CBB}
In this subsection, we fix a $\th$-twisted conjugacy class $C$ in $G$ and $w \in W$ such that $w \leq \mc$. We 
will consider the 
symplectic leaves of $\pth$ contained in $C \cap (B w B_-)$, which, by \prref{pr-prop-pi-Th}, are
the connected components of intersections of $C$ with $G^*$-orbits in $B w B_-$.
Since $\pist(h) = 0$ for every $h \in T$, it follows from 2) of \prref{pr-prop-pi-Th} that
the Poisson structure $\pTh$ on $G$ is invariant under the $\th$-twisted conjugation by elements in $T$. 
Thus for any $h \in T$, the map $h\cdot_\th: C \cap (B w B_-) \to C \cap (B w B_-)$ maps a symplectic leaf
of $\pth$ to another. If $\Sigma$ is a symplectic leaf of $\pth$, the union 
\[
T \cdot_\th \Sigma = \bigcup_{h \in T} h \cdot_\th \Sigma
\]
will be called a $T$-orbit of symplectic leaves of $\pth$, or, a {\it $T$-leaf} of $\pth$ in short.

Fix a representative $\dw \in N_G(T)$, and let $N^w = N \cap \dw N \dw^{-1}$. Then one has the unique factorization
$BwB_- = N^w T \dw N_-$ and the $(T \times T)$-action on $B w B_-$ given by
\begin{equation}\lb{eq-hhg}
(h_1, h_2) \cdot_\th (n h \dw m) = \Ad_{h_1}(n) h_1 (w\th)(h_2^{-1})h \dw 
\Ad_{\th(h_2)} (m), 
\end{equation}
where $n \in N^w, \, h \in T$, $m \in N_-$, and for $h \in T$, $(w\th)(h)= \dw(\th(h))\dw^{-1} \in T$. 
Let
\[
T_{w\th} = \{h (w\theta)(h): \;h \in T\} \hs \mbox{and} \hs 
\h_{w\th} = \{x + (w\theta)(x):\; x \in \h\}.
\]
It follows from \eqref{eq-hhg} that every $G^*$-orbit in $B w B_-$ is of the form 
\[
G^* \cdot_\theta (h_0 \dw) =N^w (\Twt h_0) \dw N_- \subset B w B_-
\]
for some 
$h_0 \in T$, and  
$\dim (G^*\cdot_\th (h_0\dw)) =\dim G - l(w) - \dim \ker (1+w\theta)$, 
where $1+w\theta: \h \to \h$, with $1$ denoting the identity operator on $\h$. 
Define
\[
\tau_{\dw}: \;\; BwB_- \lrw T, \;\; \tau_{\dw}(nh\dw m) = h, \hs n \in N^w, h \in T, m \in N_-.
\]
It follows from \eqref{eq-hhg} that 
\begin{equation}\lb{eq-taudw}
\tau_{\dw}((h_1, h_2) \cdot_\th g) = h_1 (w\th)(h_2^{-1}) \tau_{\dw}(g), \hs h_1, h_2 \in T, \, g \in B w B_-.
\end{equation}
For $h \in T$, set 
\begin{equation}\lb{eq-SwC}
\SwC = C \cap (G^* \cdot_\th (h\dw)) = \{g \in C \cap (BwB_-): \tau_{\dw}(g) \in H_{w\th}h\}.
\end{equation}
Note that $S^{\dw}_{\sC}(h_1) = S^{\dw}_{\sC}(h_2)$ if and only if $H_{w\th}h_1 = H_{w\th}h_2$.  By \eqref{eq-taudw}, one has
\begin{equation}\lb{eq-hhh}
h_1 \cdot_\th \SwC = S^{\dw}_{{\sC}}(h_1(w\th)(h_1^{-1}) h) = S^{\dw}_{\sC}(h_1^2h), \hs \forall \; h_1, h \in T.
\end{equation} 
By \prref{pr-prop-pi-Th}, symplectic leaves of $\pth$ in $C \cap (B w B_-)$ are precisely the connected components of the $\SwC$'s,
where $h \in T$.
See \cite[Remark 2.4]{E-L:grothen} for an example where $\SwC$ is not connected.

\ble{le-Gs-orbits}
1) For any $h \in T$, $\SwC \neq \emptyset$ and 
\begin{equation}\lb{eq-dim-C-Gs}
\dim \SwC= \dim C - l(w) -\dim \ker(1+w\th).
\end{equation} 

2) Let $T^\prime$ be a subtorus of $T$ such that $T^\prime T_{w\th} = T$. Then for any $h_1, h_2 \in T$ there exist $h^\prime \in T^\prime$
such that $h^\prime \cdot_\th S^{\dw}_{{\sC}}(h_1) = S^{\dw}_{{\sC}}(h_2)$;

3) For $h_0, h \in T$, $h \cdot_\th S^{\dw}_{{\sC}}(h_0) = S^{\dw}_{{\sC}}(h_0)$ if and only if 
$h^2 \in T_{w\th}$; 

4) For $h_0 \in T$, $x \in \h$ and $g \in S^{\dw}_{{\sC}}(h_0)$, $\kappa(x, x)(g) \in T_g S^{\dw}_{{\sC}}(h_0)$ if and only if
$x \in \hwt$, where $\kappa(x, x)$ is given in \eqref{eq-kappa} and $T_g \SwCz$ is the tangent space of $\SwCz$ at $g$.
\ele

\begin{proof}
We prove 2) first. Let $h_1, h_2 \in T$.
Since $T^\prime T_{w\th} = T$, there exists $h_3 \in T^\prime, h_4 \in T_{w\th}$ such that 
$h_2 = h_3h_4h_1$. Let $h^\prime \in T^\prime$ be such that $h_3 = (h^\prime)^2$. Then
$h_2 = (h^\prime)^2 h_4h_1 \in T_{w\th}(h^\prime)^2h_1$, so by \eqref{eq-hhh}, $h^\prime \cdot_\th S^{\dw}_{{\sC}}(h_1) = S^{\dw}_{{\sC}}(h_2)$.
This proves 2). Since $C \cap (B w B_-) \neq \emptyset$, 
$S^{\dw}_{{\sC}}(h_1)\neq \emptyset$ for some $h_1 \in T$. Applying 2) to $T^\prime = T$, 
one sees that $\SwC \neq \emptyset$ for all $h \in T$. 
Since the intersection $C \cap (\Gso)$ is transversal, one has \eqref{eq-dim-C-Gs}. 3) follows directly from 
\eqref{eq-hhh}. 
To prove 4), 
assume that $x \in \h$ and $g \in S^{\dw}_{{\sC}}(h_0)$ are such that $\kappa(x, x)(g) \in T_g (S^{\dw}_{{\sC}}(h_0))$.
Then there exists $y^* = (y_+ + y, \, -y+y_-) \in \g_{\rm st}^*$, where $y_+ \in \n, y_- \in \n_-$ and $y \in \h$, such that $y^*-(x, x) 
\in \Ad_{(g, e)} \{(z, \theta^{-1}(z)): z \in \g\}$, or
$y_+ + y -x = \Ad_g \theta(y_- - y - x)$. Write $g = n hh_0 \dw n_-$ for some $n \in N^w, h \in \h_{w\th}$, and
$n_- \in N_-$. One has
\[
\Ad_{nhh_0}^{-1}(y_+ + y-x) = \Ad_{\dw} \Ad_{n_-} \theta(y_- - y - x) \in \b \cap \Ad_{\dw} \b_-=\h + \n \cap \Ad_{\dw} \n_-,
\]
whose $\h$-component gives $y-x = -w\theta(y+x)$, and thus
\[
x = \frac{1}{2} (x-w\theta(x)) + \frac{1}{2}(x + w\theta(x)) = \frac{1}{2}(x+y + w\theta(x+y)) \in \hwt.
\]
\end{proof}

\bpr{pr-one-T-leaf} Let $T^\prime$ be a subtorus of $T$ such that $T^\prime T_{w\th} = T$. Then 
 $C \cap (B w B_-)$ is one single $T^\prime$-orbit of symplectic leaves of $\pth$ under the $\th$-twisted conjugation.
\epr

\begin{proof} Consider the action map
\[
\sigma: \;\; T^\prime \times \Sigma \longrightarrow C \cap (B w B_-), \;\; (h, g) \longmapsto h g \th(h^{-1}), \hs h \in T^\prime, \, g \in 
\Sigma.
\]
Let $\h^\prime \subset \h$ be the Lie algebra of $T^\prime$. For any $g \in \Sigma$,
the kernel of the differential $\sigma_*(e, g)$ of $\sigma$ at $(e, g)$ is isomorphic to
$\{x \in \h^\prime: \kappa(x, x)(g) \in T_g \Sigma\}$, which, by \leref{le-Gs-orbits}, is $\h' \cap \h_{w\th}$.
Let $\h^{\prime\prime}$ be any complement of $\h^\prime \cap \h_{w\th}$ in $\h^\prime$. It follows from 
$T^\prime T_{w\th} = T$ that $\h^\prime + \h_{w\th} = \h$ and thus $\h = \h^{\prime\prime} + \h_{w\th}$ is a direct sum.
Thus $\sigma_*(e, g)$ maps $\h^{\prime\prime} \times T_g \Sigma$ injectively to $T_g(C \cap (B w B_-))$. Since
 $\dim (\h^{\prime\prime} \times T_g \Sigma) = \dim T_g(C \cap (B w B_-))$, $\sigma$ is a submersion at $(e, g)$. Since
$T$ is abelian, $\sigma$ is a submersion everywhere. Thus $T^\prime \cdot_\th \Sigma$
is an open subset of
$C \cap (B w B_-)$. If $\Sigma^\prime$ is another symplectic leaf of $\pth$ in $C \cap (B w B_-)$, 
$T^\prime \cdot_\th \Sigma^\prime$
is also an open subset of $C \cap (B w B_-)$, and either $T^\prime \cdot_\th \Sigma =T^\prime \cdot_\th \Sigma^\prime$
or
$(T^\prime \cdot_\th \Sigma) \cap (T^\prime \cdot_\th \Sigma^\prime) = \emptyset$. Since $C \cap (B w B_-)$ is connected, one must have 
$T^\prime \cdot_\th \Sigma =T^\prime \cdot_\th \Sigma^\prime$. This shows that $T^\prime \cdot_\th \Sigma  
= C \cap (B w B_-)$.
\end{proof}

\bco{co-leaves}
With the notation as in \prref{pr-one-T-leaf}, let $h \in T$ and let $\Sigma$ be any connected component of
$\SwC$. Then $\Sigma$ is a symplectic leaf of $\pth$ in $C \cap (B w B_-)$, and every symplectic leaf of
$\pth$ in $C \cap (B w B_-)$ is of the form $h'\cdot_\th \Sigma$ for some $h^\prime \in T^\prime$.
\eco

\subsection{The $T$-leaves of $\pTh$}\lb{subsec-T-orb}

Summarizing the results in $\S$\ref{subsec-CBB} and $\S$\ref{subsec-leaves-CBB}, we have the following \thref{th-TT-orb}.

\bth{th-TT-orb}  Let $T^\prime$ be a subtorus of $T$ such that $T^\prime T_{w\th} = T$ for every $w \in W$.
Then the $T^\prime$-orbits of symplectic leaves of $\pth$ in $G$ are precisely all the intersections
$C \cap (BwB_-)$, where $C$ is a $\th$-twisted conjugacy class in $G$
and $w \in W$ is such that $w \leq \mc$. Every such an intersection $C \cap (BwB_-)$ is non-empty and the symplectic leaves in  
$C \cap (BwB_-)$ have dimension  equal to $\dim C - l(w) -\dim \ker (1+w\theta)$.
\eth

\thref{th-T-orb-mc} is now a special case of \thref{th-TT-orb} by taking $T^\prime = T$.

In $\S$\ref{sec-ptth}, we apply \thref{th-TT-orb} to an example where  $T^\prime$ is a proper subtorus of  $T$.

\sectionnew{The lowest rank of $\pth$ in a $\th$-twisted conjugacy class}\lb{sec-lowest-C}

\subsection{Some linear algebra facts}\lb{subsec-linear-algebra}
Recall that for an $n$-dimensional real or complex vector space $V$, $S \in {\rm End}(V)$ is called a 
reflection if $\dk(1-S) = n-1$ and $\dk(1+S) = 1$,
where $1$ also
denotes the identity map on $V$.

\ble{le-linear}
Let $V$ be an $n$-dimensional real vector space and let $S \in {\rm End}(V)$ be
a reflection.  Then for any 
$A \in {\rm End}(V)$, 
\begin{align}\lb{eq-rk-dk-1}
-1 + \dim \ker  (1+A)& \leq \dim \ker  (1+AS) \leq 1+ \dim \ker  (1+A);\\
\lb{eq-rk-dk-2}
-1 + \rk (1+A)& \leq \rk(1+AS) \leq 1+ \rk(1+A).
\end{align}
Moreover,  $\rk(1+A)=1+ \rk(1+AS)$ if and only if $\ker(1+S) \subset {\rm im}(1+A)$ and 
$\ker(1+S) \cap {\rm im}(1+AS)=0$.
\ele

\begin{proof} Let $x_0  \in \ker(1+S)$, $x_0 \neq 0$. Then ${\rm im}(1-S) = {\mathbb R}x_0$. It follows from
$1+AS = (1+A) -A(1-S)$ that ${\rm im}(1+AS) \subset {\rm im}(1+A) + {\mathbb R} x_0$. Exchanging $A$ and $AS$, one has
\begin{equation}\lb{eq-AAS}
{\rm im} (1+A) + {\mathbb R} x_0 = {\rm im}(1+AS) + {\mathbb R} x_0,
\end{equation}
from which both inequalities in  
\eqref{eq-rk-dk-2} follow.
Using $\dim \ker  (A) = n-\rk (A)$ for  $A \in {\rm End}(V)$, one proves \eqref{eq-rk-dk-1}. The last statement 
also follows from \eqref{eq-AAS}.
\end{proof}

\ble{le-linear-1}
Let $V$ be a finite dimensional real vector space with an inner product $\lara$, and let $O(V, \lara)$ be
the group of all linear operators on $V$ preserving $\lara$. 
Then for all $A, S \in O(V, \lara)$ and $S$ a reflection, one has 
\[
\dk (1+AS) = \dk (1+A)\pm 1 \hs \mbox{and} \hs \rk (1+AS) = \rk (1+A)\pm 1.
\]
\ele

\begin{proof} If $\lambda \in \Cset$ is an eigenvalue for
$A$, so is $\bar{\lambda}$, and $|\lambda|^2 = 1$. It follows that $\det(A) = (-1)^{\dim \ker (1+A)}$. Since 
$\det(AS) = -\det(A)$, $\dim \ker (1+AS) \neq \dim \ker (1+A)$. It follows from \leref{le-linear} that 
$\dim \ker (1+AS) = \dim \ker (1+A) \pm 1$, and that 
\[
\rk (1+AS) = \dim V - \dk(1+AS) = \rk (1+A)\pm 1.
\]
\end{proof} 

\ble{le-inv}
For any finite dimensional real vector space $V$ and $A \in {\rm End}(V)$, one has
$\ker (1+A) \subset \im (1-A)$ and $\rk(1-A)-\dim \ker  (1+A) = \rk(1-A^2)$.
\ele

\begin{proof}
Since $2x = (x+Ax) + (x-Ax)$ for every $x \in V$, one has $\ker (1+A) \subset \im (1-A)$. Since the
map $(1+A)|_{{\rm im}(1-A)}: {\rm im}(1-A) \to {\rm im}(1-A^2)$ is surjective and has $\ker (1+A)$ as its kernel, one
has $\rk(1-A)-\dim \ker  (1+A) = \rk(1-A^2)$.
\end{proof}

We now apply the above linear algebra facts to 
the case of $V = \h$ with the inner product $\lara$ induced from the Killing form of $\g$.  
We continue to assume that $\th(B) = B$ and $\th(T) = T$, so 
$\th \in O(\h, \lara)$. For $w \in W$, let 
\[
L_\th(w) = l(w) + \dim \ker (1+w\theta) \hs \mbox{and} \hs L_\th^\prime(w) = l(w) + \rk(1-w\theta).
\]
By \leref{le-inv}, one has
\begin{equation}\lb{eq-LL-ww}
L_\th^\prime(w) - L_\th(w) = \rk(1-(w\th)^2), \hs w \in W.
\end{equation}
For a root $\alpha \in \h^*$, let $s_\al \in W$ be the corresponding reflection on $\h$.

\ble{le-wwp}
If $w, w' \in W$ are such that $w \leq w'$, then
\begin{equation}\lb{eq-rk-wwp}
L_\th(w) \leq L_\th(w')  \hs \mbox{and} \hs 
L_\th^\prime(w) \leq L_\th^\prime(w'),
\end{equation}
and both $L_\th(w')-L_\th(w)$ and $L_\th^\prime(w')-L_\th^\prime(w)$ are even integers.
\ele

\begin{proof} 
By \cite[Corollary 2.2.2]{B:flag}, there exists a sequence $(w_1, \ldots, w_k)$ in $W$
such that $w = w_1 < \cdots < w_k = w'$ and $l(w_{j+1}) = l(w_j) + 1$ for $1 \leq j \leq k-1$. We can thus assume that $l(w') = l(w) + 1$. 
Then $w' = ws_\alpha$ for a positive root $\alpha$. By \leref{le-linear-1},
\[
L_\th(w) =l(w) + \dk(1+w'\theta s_{\theta^{-1}(\alpha)})
= l(w') -1+\dk(1+w'\theta) \pm 1,
\]
so $L_\th(w')-L_\th(w)\in \{0, 2\}$. 
Similarly, one has  $L_\th^\prime(w')-L_\th^\prime(w) \in \{0, 2\}$.
\end{proof}

\subsection{The lowest rank of $\pth$ in a $\th$-twisted conjugacy class}
For $g \in G$, denote by $\rk(\pTh(g))$ the dimension of
the symplectic leaf of $\pTh$ through $g$. By \thref{th-T-orb-mc}, 
\begin{equation}\lb{eq-rk-g}
\rk(\pth(g)) = \dim C - L_\th(w), \hs \forall \; g \in C \cap (B w B_-),
\end{equation} 
where $C$ is any $\th$-twisted conjugacy class
and $w \in W$ is such that $w \leq \mc$. 
By \leref{le-wwp}, if $w, w' \in W$ are such that $w \leq w'$, then 
\begin{equation}\lb{eq-rank-g-g}
\rk(\pth(g)) \geq \rk(\pth(g')), \hs \mbox{if} \;\; \forall \; g \in C \cap (BwB_-), \,\, g' \in C \cap (B w' B_-).
\end{equation}

\bpr{pr-min-rank-C} For any $\th$-twisted conjugacy class $C$, 
\[
{\rm min}\{\rk(\pTh(g)): g \in C\} =\dim C - L_\th(\mc),
\]
and for any $g \in C$,  $\rk(\pth(g)) = \dim C - L_\th(\mc)$ if and only if  
$g \in C \cap (B\mc B_-)$.
\epr

\begin{proof} 
By \eqref{eq-rank-g-g}, $\dim C - L_\th(\mc)=\rk(\pth(g))$ for any $g \in C \cap (B \mc B_-)$ is the minimal rank of $\pth$ in $C$. 
It remains to show that if $w \in W$ is such that $w < \mc$, then $L_\th(w) \leq  L_\th(\mc)-2$. 

Let $w \in W$ be such that
$w < \mc$. By \cite[Corollary 2.2.2]{B:flag}, there exists a sequence $(w_1, \ldots, w_k)$ in $W$
such that $w = w_1 < \cdots < w_k = \mc$ and $l(w_{j+1}) = l(w_j) + 1$ for $1 \leq j \leq k-1$. By \eqref{eq-rank-g-g}, 
$\rk(\pth(g)) \geq \dim C - L_\th(w_{k-1})$ for all $g \in C \cap (B w B_-)$. It thus suffices to show that $L_\th(w) = L_\th(\mc)-2$ for any $w \in W$
such that $w < \mc$ and $l(w) = l(\mc)-1$.

Assume thus that $w < \mc$ and $l(w) = l(\mc)-1$. Then $w = \mc s_\beta$ for some positive root $\beta$ such that
$\mc(\beta) < 0$. By the proof of \leref{le-wwp}, $L_\th(w) = L_\th(\mc)$ or $L_\th(w) = L_\th(\mc)-2$. Suppose that
$L_\th(w) = L_\th(\mc)$. Then 
\[
\rk(1+\mc \th) = 1 + \rk(1+w\th) = 1 + \rk(1+\mc \th s_{\th^{-1}(\beta)}).
\]
Regarding both $\mc$ and $\th$ as linear operators on $\h^*$, one sees from \leref{le-linear} that $\th^{-1}(\beta) \in {\rm im}(1+\mc \th)$,
and thus by \eqref{eq-mc-inv},
\[
\th^{-1}(\beta)-\mc(\beta) =(1-\mc \th) \th^{-1}(\beta) \in {\rm im} (1-(\mc \th)^2)= {\rm im}(1-\th^2).
\]
Let $\lambda \in \h^*$ be
such that $\th^{-1}(\beta) - \mc(\beta) = (1-\th^2) (\lambda).$
Let $\{\al_1, \ldots, \al_k\}$ be the set of all simple roots and write $\lambda = \sum_{j = 1}^k a_j \al_j$ and 
$\th^2(\lambda) = \sum_{j=1}^k b_j \al_j$. Then 
\begin{equation}\lb{eq-lam-beta}
\th^{-1}(\beta) - \mc(\beta) = \sum_{j=1}^k (a_j-b_j)\al_j.
\end{equation}
Since both $\th^{-1}(\beta)$ and $-\mc(\beta)$ are positive roots, each $a_j - b_j$ is a non-negative integer and 
$\sum_{j=1}^k (a_j-b_j) > 0$. On the other hand, since $\th^2$ permutes the simple roots, 
$\sum_{j=1}^k a_j =\sum_{j=1}^k b_j$. We have thus arrived at a contradiction and we conclude that $L_\th(w) = L_\th(\mc)-2$.
\end{proof}

\bre{re-min-rank} From \prref{pr-min-rank-C} and  
the linear algebra fact that if $t \to A(t) \in {\rm End}(V)$ is a continuous map from an open 
neighborhood of $0 \in {\mathbb C}$ to ${\rm End}(V)$ then $\rk(A(t)) \geq \rk(A(0))$ for $t$ closed to $0$, one concludes
that $C \cap (B \mc B_-)$ is closed in $C$. To see this from another way, observe first that  
\[
B \mc B_- \subset \bigsqcup_{w \in W, w\geq \mc} BwB
\]
 by \cite[Corollary 1.2]{De}, so by the
definition of $\mc$, 
$C \cap (B\mc B_-) \subset B \mc B.$
Since 
\[
(B \mc B_-) \cap (B \mc B) = B \mc,
\]
one has
$C \cap (B \mc B_-) = C \cap (B \mc)$. 
In particular, $C \cap (B\mc B_-)$ is closed in $C$. Similarly, $C \cap (B \mc B_-) = C \cap (\mc B_-)$. 
\ere

\sectionnew{Spherical conjugacy classes and the lowest rank of $\pth$ in $G$}\lb{sec-spherical}

\subsection{Spherical $\theta$-twisted conjugacy classes} 
Recall
that a $\theta$-twisted conjugacy class $C$ in $G$ is said to be spherical if it has an open orbit
under the $\theta$-twisted conjugation by some Borel subgroup of $G$. 
For $g \in G$, let $B \cTh g$ be the 
$B$-orbit through $g$ under the $\th$-twisted conjugation.
The following \leref{le-stabi-B} on the dimensions of $B$-orbits in $G$ is a generalization of \cite[Theorem 5]{CCC}
and is proved in \cite{Co:bad}.

\ble{le-stabi-B}
For any $w \in W$ and $g \in BwB$, one has 
\[
\dim (B \cTh g) \geq l(w) + \rk (1-w\theta).
\]
In particular, $\dim C \geq l(\mc) + \rk(1-\mc\theta)$ for any $\theta$-twisted conjugacy class $C$.
\ele

Generalizing a result in \cite{CCC} of 
N. Cantarini, G. Carnovale, and M. Costantini for $\theta={\rm Id}_G$, the following characterization of
spherical $\theta$-twisted 
conjugacy classes is proved in \cite{Lu:formula} for the case of characteristic zero and in \cite{Car-2012} 
for good odd characteristics.

\bth{th-dim-formula}\cite{CCC, Car-2012, Lu:formula} For $\theta \in \Aut(G)$ such that $\th(B) = B$ and $\th(T) = T$,
a $\theta$-twisted conjugacy class $C$ in $G$ is spherical if and only if
\[
\dim C = l(\mc) + \rk(1-\mc\theta).
\]
\eth

\bco{co-spherical-min-rank}
If $C$ is a spherical $\th$-twisted conjugacy class in $G$, then
\[
{\rm min}\{\rk(\pth(g)): g \in C\} = \rk(1-\th^2),
\]
and for $g \in C$, $\rk(\pth(g)) = \rk(1-\th^2)$ if and only if $g \in C \cap (B \mc B_-)$.
\eco

\begin{proof} If $C$ is spherical, then 
\[
\dim C - L_\th(\mc) = L_\th^\prime(\mc) - L_\th(\mc) = \rk(1-(\mc \th)^2) = \rk(1-\th^2),
\]
and \coref{co-spherical-min-rank} now follows directly from \prref{pr-min-rank-C}.
\end{proof}

\bre{re-rk-even}
Note that as $\th$ is a linear operator on $\h$ preserving the inner product on $\h$ induced by the Killing form
of $\g$, one can decompose $\h$ into the direct sum of $\ker(1-\th^2)$ and $2$-dimensional subspaces on which 
$\th^2$ acts as non-trivial rotations. Consequently, $\rk(1-\th^2)$ is always an even integer.
\ere

\bex{ex-triality}
Assume that $G$ is of type  $D_4$ and $\th \in {\rm Aut}(G)$ has order $3$ such that
the induced $\th \in {\rm Aut}(\Gamma)$ fixes $\al_2$ and maps $\al_1$ to $\al_3$, $\al_3$ to $\al_4$, 
and $\al_4$ to $\al_1$, where $\al_2$ is the simple root not orthogonal to any of the other three. Then the fixed point set $G^\th$ of $\th$
is of type $G_2$ and the $\th$-twisted conjugacy class $C$ through the identity element of $e$ is spherical \cite[$\S$4.5]{Co:bad}. 
In this case (see \cite[$\S$4.5]{Co:bad} and \cite[Example 3.9]{Lu:formula}), $\mc = w_0s_{\al_2}$, where $w_0$ is the longest element in $W$, and one has
$\dim C =14 = l(w_0s_{\al_2})+\rk(1-w_0s_{\al_2}\th)$. Since $\rk(1-\th^2) = 2$, the lowest rank of $\pth$ in $C$ is
$2$.
\eex

\bre{re-spherical-exist-1}
When $\th^2 = 1 \in {\rm Aut}(\Gamma)$, or when $G$ is simple, spherical $\th$-twisted conjugacy classes exist for any $\th \in \Aut(G)$.
Indeed, by \reref{re-Th-Th}, we may assume that $\th(B) = B$ and $\th(T) = T$.
Fix a root vector $e_\alpha$ for each $\alpha \in 
\Gamma$ and let $\lambda \in {\mathbb C}$ be such that $\theta(e_\alpha) = \lambda_\alpha e_{\theta(\alpha)}$.
Choose any $h \in T$ such that $h^\alpha = \frac{1}{\lambda_\alpha}$ for each $\alpha \in \Gamma$, and let 
$\theta_1 =\th \Ad_h =\Ad_{\th(h)} \th\in \Aut(G)$. Then $\th_1(\al) = \th(\al)$ and $\theta_1 (e_\alpha) = e_{\theta(\alpha)}$
for every $\alpha \in \Gamma$, so $\th_1 \in {\rm Aut}(G)$ has the same order as $\th \in {\rm Aut}(\Gamma)$.
Suppose that $\th^2 = 1 \in {\rm Aut}(\Gamma)$. Then 
$\theta_1^2 = {\rm Id}_G$,  and the $\theta_1$-twisted conjugacy class through the identity element is trivial if $\th_1 = {\rm Id}_G$ and 
a symmetric space, and hence spherical, if $\th_1$ has order $2$, and thus the $\theta$-twisted conjugacy class $C$ through $\th(h)$ is spherical.
If $G$ is simple, $\th \in {\rm Aut}(\Gamma)$ has order $2$ or $3$, with the latter 
(triality case) possible only for 
$G$ of type $D_4$, and we saw in \exref{ex-triality} that spherical $\th$-twisted conjugacy classes exist in this case.  
When $G$ is simple, spherical conjugacy classes in $G$ (i.e. the case for $\th = {\rm Id}_G$) 
have been classified 
in \cite{CCC, Panyushev}. 
\ere

\subsection{Proof of \thref{th-min-rank}} 
For any $\th$-twisted conjugacy class $C$ in $G$ and any $g \in C$, one has 
\[
\rk(\pth(g)) \geq \dim C - L_\th(\mc) \geq L_\th^\prime(\mc) - L_\th(\mc) =\rk(1-\th^2).
\]
We have already proved in \coref{co-spherical-min-rank} that $\rk(\pth(g)) = \rk(1-\th^2)$ for all $g \in C \cap (B \mc B_-)$
if $C$ is spherical.
Suppose that $g \in G$ is such that $\rk(\pth(g)) = \rk(1-\th^2)$ and let $C$ be the $\th$-twisted conjugacy class
of $g$. By \prref{pr-min-rank-C}, $\rk(\pth(g_1)) = \rk(1-\th^2)$ for all $g_1 \in C \cap (B \mc B_-)$, and thus 
\[
\dim C = L_\th(\mc) + \rk(1-\th^2) = L_\th^\prime(\mc).
\]
By \thref{th-dim-formula}, $C$ is spherical, and by \prref{pr-min-rank-C}, $g \in C \cap (B \mc B_-)$.
\qed

\subsection{Proof of \coref{co-0-locus}} \coref{co-0-locus} follows directly from \thref{th-T-orb-mc}, \thref{th-min-rank} and the statement in
\reref{re-spherical-exist-1} that spherical $\th$-twisted conjugacy classes exist when $\th^2 =1 \in {\rm Aut}(\Gamma)$.

\bex{ex-SL-2}
For $G = SL(n+1, {\mathbb C})$, a conjugacy class $C$ in $G$ is spherical if and only if either 
$C = C^{\rm ss}_{l, \lambda, \lambda'}$, where $1 \leq l \leq (n+1)/2$ is an integer and 
$C^{\rm ss}_{l, \lambda, \lambda'}$ is the semisimple conjugacy class with two distinguished eigenvalues $\lambda$ and $\lambda^\prime$
of respective multiplicities $l$ and $n+1-l$, or $C = C^{\rm uni}_{l, \lambda}$, where $0 \leq l \leq (n+1)/2$ is an integer, and
$C^{\rm uni}_{l, \lambda}$ is the conjugacy class of $\lambda g$ with $g$ an 
upper triangular Jordan matrix having $1$ as the only eigenvalue and  $n+1-l$ size $1$ Jordan blocks and $l$ size $2$ Jordan blocks
(see \cite{Can, CCC}). 
Let again $B$ and $B_-$ be the subgroups of $G$ consisting of all upper-triangular and lower triangular matrices respectively. For 
$C= C^{\rm ss}_{l, \lambda, \lambda'}$ or $C = C^{\rm uni}_{l, \lambda}$ as above, it follows from \exref{ex-SL} that 
$\mc$ has the cycle decomposition
\[
\mc = m_l =(1, n+1) (2, n) \cdots (l,\; n+2-l),
\]
and thus $\dim (C \cap B \mc B_-) = \dim C - l(\mc) = \rk(1-\mc) = l$, so by \coref{co-0-locus}, $C \cap (B \mc B_-) \cong ({\mathbb C}^\times)^l$.
Explicit parameterizations of  $C \cap (B \mc B_-)$ by $({\mathbb C}^\times)^l$ were given in \cite{Simon-thesis}, where the 
Poisson structure $\pth$ for $G = SL(n+1, {\mathbb C})$ and $\th = {\rm Id}_G$ was treated in more detail. For an integer $k$, let $I_k$ be
the $k \times k$ identity matrix and $J_k$ the $k \times k$ anti-diagonal matrix with $1$'s on the anti-diagonal. Let
$0 \leq l \leq (n+1)/2$ be an integer and let 
$\lambda, \lambda^\prime \in {\mathbb C}$ be such that $\lambda^{n+1-l} (\lambda^\prime)^l = 1$.  
For $(x_1, \ldots, x_l) \in ({\mathbb C}^\times)^l$, let 
\[
g(\lambda, \lambda^\prime, x_1, \ldots, x_l) = \left(\begin{array}{ccc} (\lambda+\lambda^\prime) I_l & 0 & \lambda^\prime XJ_l\\
0 & \lambda I_{n+1-2l} & 0 \\ -\lambda X^\prime J_l & 0 & 0\end{array}\right) \in SL(n+1, {\mathbb C}),
\]
where $X = {\rm diag}(x_1, x_2, \ldots, x_l)$ and $X^\prime = {\rm diag}(x_l^{-1}, \ldots, x_2^{-1}, x_1^{-1})$.
By \cite[$\S$5.3]{Simon-thesis},
\begin{align*}
(C^{\rm ss}_{l, \lambda, \lambda'}) \cap (B \mc B_-) &=\{g(\lambda, \lambda^\prime, x_1, \ldots, x_l): \,(x_1, \ldots, x_l) \in ({\mathbb C}^\times)^l
\},\\
(C^{\rm uni}_{l, \lambda}) \cap (B \mc B_-) & =  \{g(\lambda, \lambda, x_1, \ldots, x_l): \,(x_1, \ldots, x_l) \in ({\mathbb C}^\times)^l
\}.
\end{align*}
We conclude that the zero-locus $Z(\pth)$ of $\pth$ in $SL(n+1, {\mathbb C})$ is given by
\[
Z(\pth) = \bigcup_{l=0}^{[(n+1)/2]} \{g(\lambda, \lambda^\prime, x_1, \ldots, x_l): \,\lambda^{n+1-l} (\lambda^\prime)^l = 1, \, (x_1, \ldots, x_l) \in ({\mathbb C}^\times)^l\}.
\]
\eex

\sectionnew{The Poisson structure $\pth$ for an arbitrary $\th \in {\rm Aut}(G)$}\lb{sec-th-arbitrary}
Recall that we fixed, from the very beginning, the pair $(B, T)$ of a Borel subgroup and a maximal torus $T \subset B$.
Although the Poisson structure $\pth$ on $G$ is defined in $\S$\ref{subsec-dfn} for any $\th \in {\rm Aut}(G)$ and \prref{pr-prop-pi-Th}
holds for $\pth$, we made the assumption that $\th$ stabilizes both $B$ and $T$ throughout $\S$\ref{sec-T-leaves}-$\S$\ref{sec-spherical}.

Assume now that $\th \in {\rm Aut}(G)$ is arbitrary. By \reref{re-Th-Th}, there exists $g_0 \in G$ such that $\th(B) = g_0 B g_0^{-1}$
and $\th(T) = g_0 T g_0^{-1}$. Choose any such $g_0$ and let $\th'=\Ad_{g_0^{-1}} \th$. Then $\th'$
stabilizes both $B$ and $T$, and 
\[
r_{g_0}: \;\; (G, \pth) \longrightarrow (G, \pi_{\th^\prime}), \;\; g \longmapsto gg_0, \hs g \in G,
\]
is a Poisson isomorphism. Applying results in $\S$\ref{sec-T-leaves}-$\S$\ref{sec-spherical} to $\pi_{\th'}$, we 
arrive at properties of $\pth$ via the Poisson isomorphism $r_{g_0}$. In particular, \thref{th-T-orb-mc}, \thref{th-min-rank}, 
and \coref{co-0-locus} hold if $B w B_-$ is replaced by $B w B_- g_0^{-1}$ for $w \in W$.

\sectionnew{An example: the Poisson structure $\ptth$ on $G^n$}\lb{sec-ptth}

\subsection{The Poisson structure $\ptth$ on $G^n$}
In this section, we assume that $G$ is non-trivial and consider the direct product group $\wt{G} = G^n$, where $n \geq 2$, and 
$\tth\in 
{\rm Aut}(G^n)$ is given by 
\[
\tth(g_1, g_2, \ldots, g_n) = (g_2, g_3, \ldots, g_n, g_1), \hs g_j \in G, \, j = 1, 2, \ldots, n.
\]
Let  $\widetilde{B} = B^n$, $\widetilde{B}_- = B_-^n$, and $\widetilde{T} = T^n$.
Then $\tth(\widetilde{B}) = \widetilde{B}$ and $\tth(\widetilde{T}) = \widetilde{T}$. Applying the construction 
in $\S$\ref{subsec-dfn} to $G^n$ using the pair $(\widetilde{B}, \wt{T})$ and $\tth \in {\rm Aut}(G^n)$, one obtains the holomorphic Poisson structure
$\ptth$ on $G^n$. We describe some properties of $\ptth$.

The $\tth$-twisted conjugation of $G^n$ on itself is given by
\begin{equation}\lb{eq-gg-kk}
(g_1, g_2, \ldots, g_n) \cdot_{\tth} (k_1, k_2, \ldots, k_n) = (g_1k_1g_2^{-1}, \, g_2k_2 g_3^{-1}, \, \ldots, g_nk_ng_1^{-1}).
\end{equation}
Consider the multiplication map
\[
\mu_n:\;\;\;  G^n \longrightarrow G,\;\;\; (g_1, g_2, \ldots, g_n) \longmapsto g_1g_2 \cdots g_n.
\]
For a conjugacy class $C$ in $G$, let $\widetilde{C} = \mu_n^{-1}(C) \subset G^n$. 

\ble{le-mun-C}
$\tth$-twisted conjugacy classes
in $G^n$ are precisely the subsets $\wt{C}$ of $G^n$, where $C$ is a conjugacy class in $G$. 
\ele

\begin{proof}
Let $\tilde{k} = (k_1, k_2, \ldots, k_n) \in G^n$, let ${\mathcal O}_{(k_1,k_2, \ldots,k_n)}$ be the $\tth$-twisted 
conjugacy class of $\tilde{k}$, and let $C_{k_1k_2 \cdots k_n}$ be the conjugacy class of $k_1k_2 \cdots k_n$ in $G$. It is clear from 
\eqref{eq-gg-kk} that ${\mathcal O}_{(k_1,k_2, \ldots,k_n)} \subset \mu_n^{-1}(C_{k_1k_2 \cdots k_n})$. Conversely, 
let $\tilde{h} = (h_1, h_2, \ldots, h_n) \in \mu_n^{-1}(C_{k_1k_2 \cdots k_n})$
and let $g_1 \in G$ be such that $h_1h_2 \cdots h_n = g_1k_1k_2 \cdots k_n g_1^{-1}$. For $2 \leq j \leq n$, let
$g_j \in G$ be given inductively by $g_{j-1}k_{j-1}g_j^{-1} = h_{j-1}$. It follows from
$h_1h_2 \cdots h_n = g_1k_1k_2 \cdots k_n g_1^{-1}$ that $h_n = g_n k_n g_1^{-1}$, and thus
$\tilde{g} \cdot_{\tth} \tilde{k} = \tilde{h}$, where $\tilde{g} = (g_1, g_2, \ldots, g_n)$. This shows that
$\mu_n^{-1}(C_{k_1k_2 \cdots k_n}) \subset {\mathcal O}_{(k_1,k_2, \ldots,k_n)}$. Hence 
$\mu_n^{-1}(C_{k_1k_2 \cdots k_n})
= {\mathcal O}_{(k_1,k_2, \ldots,k_n)}$, 
and \leref{le-mun-C} follows.
\end{proof}

Identify the Weyl group $\wt{W}$ of $G^n$ with $W^n$. Then for $\tilde{w} = (w_1, w_2, \ldots, w_n) \in W^n$,
\[
\widetilde{B} \tilde{w} \widetilde{B}_-= (B w_1B_-) \times (B w_2 B_-) \times \cdots \times (B w_n B_-).
\]
Recall that we assume that $n \geq 2$.

\ble{le-tilde-W}
For any conjugacy class $C$ in $G$ and any $\tilde{w} = (w_1, w_2, \ldots, w_n) \in W^n$, 
\[
\widetilde{C} \cap (\widetilde{B} \tilde{w} \widetilde{B}_-) \neq \emptyset.
\]
\ele

\begin{proof}
Let $C$ be a conjugacy class in $G$.
By \prref{pr-mc}, it is enough to show that $\widetilde{C} \cap (\widetilde{B} \tilde{w}_0 \widetilde{B}_-) \neq \emptyset$, where
$\tilde{w}_0 = (w_0, w_0, \ldots, w_0) \in W^n$, and $w_0$ is the longest element in $W$. Since $\widetilde{B} \tilde{w}_0 \widetilde{B}_-
=(Bw_0)^n$, one needs to show that $C \cap \mu_n((Bw_0)^n) \neq \emptyset$.

If $n \geq 2$ is even, since $Bw_0Bw_0 = BB_- \ni e$, one has
\[
\mu_n((Bw_0)^n) \supset Bw_0Bw_0 = BB_- \supset B,
\]
and since $C \cap B \neq \emptyset$, one has $C \cap \mu_n((Bw_0)^n) \neq \emptyset$. Suppose that $n \geq 3$ is odd. 
Then $\mu_n((Bw_0)^n) \supset Bw_0Bw_0Bw_0 = Bw_0BB_-$. We claim that $Bw_0BB_-=G$ and thus $C \cap \mu_n((Bw_0)^n) \neq \emptyset$. 
To show that $Bw_0BB_-=G$, first note that $Bw_0BB_-$ is the union of some
$(B, B_-)$-double cosets. For any $w \in W$, since $w \leq w_0$, $BwB_- \cap Bw_0B \neq \emptyset$ by \cite[Corollary 1.2]{De},
so $B w B_- \subset Bw_0BB_-$. Since $w \in W$ is arbitrary, one has $Bw_0BB_-=G$.
\end{proof}
 
In the notation of $\S$\ref{subsec-T-orb}, for every conjugacy class $C$ in $G$ one has 
$\wt{W}_{\scriptscriptstyle \widetilde{C}}^- = \wt{W}$ and $m_{\widetilde{\scriptscriptstyle C}} = (w_0, w_0, \ldots, w_0)$.
Applying \thref{th-T-orb-mc}, we have the following description of $\wt{T}$-orbits of symplectic leaves
of $\ptth$ on $G^n$, where $\wt{T}$ acts on $G^n$ by the $\tth$-twisted conjugation given in \eqref{eq-gg-kk}.

\bpr{pr-T-orb-ptth}
For any conjugacy class $C$ and $\tilde{w} = (w_1, w_2, \ldots, w_n) \in W^n$, 
\begin{align*}
\wt{C} \cap (\widetilde{B} \tilde{w} \widetilde{B}_-) = &\{(g_1, g_2, \ldots, g_n) \in 
(B w_1B_-) \times (B w_2 B_-) \times \cdots \times (B w_n B_-): \\
& \;\;g_1g_2 \cdots g_n \in C\}
\end{align*}
is nonempty, and  the $\wt{T}$-orbits of symplectic leaves of $\ptth$ in $G^n$ are precisely all such intersections;
The dimension of the symplectic leaves of $\ptth$ in 
$\wt{C} \cap (\widetilde{B} \tilde{w} \widetilde{B}_-)$ is 
\[
\dim C + (n-1)\dim G -(l(w_1) + l(w_2) + \cdots + l(w_n)) - \dim \ker (1-(-1)^n w_1w_2 \cdots w_n),
\] 
where $1-(-1)^n w_1w_2 \cdots w_n$ is regarded as a linear operator on the Lie algebra $\h$ of $T$.
\epr

\begin{proof}
We only need to prove the last statement on the dimension of the symplectic leaves in 
$\wt{C} \cap (\widetilde{B} \tilde{w} \widetilde{B}_-)$, which by \thref{th-T-orb-mc}, is equal to 
\[
\dim \wt{C} -(l(w_1) + l(w_2) + \cdots + l(w_n))-\dim \ker (1+(w_1, w_2, \ldots, w_n)\tth),
\]
where $1+(w_1, w_2, \ldots, w_n)\tth$ is regarded as a linear operator on the Lie algebra of $\wt{T}$. It is easy to see
that $\ker (1+(w_1, w_2, \ldots, w_n)\tth)\cong \ker (1-(-1)^n w_1w_2 \cdots w_n)$, and since 
$\wt{C} \cong G^{n-1} \times C$, so $\dim \wt{C} = \dim C + (n-1)\dim G$, and the claim is proved.
\end{proof}

\subsection{Spherical $\tth$-twisted conjugacy classes}
We now classify spherical $\tth$-twisted conjugacy classes in $G^n$. Recall that we assume $G$ to be non-trivial.

\ble{le-tth-spherical}
1) If $n >3$ or if $n = 3$ and the rank of $G$ is greater than $1$, there are no spherical $\tth$-twisted conjugacy classes in $G^n$;

2) For $n = 3$ and $G$ having rank $1$, or for $ n= 2$ and $G$ arbitrary, $\wt{C}$ is spherical if and only if 
$C =\{z\}$, where $z$ is in the center of $G$.
\ele

\begin{proof}
We already know that   
$\dim \wt{C} = \dim C + (n-1) \dim G$ for any conjugacy class $C$ in $G$. 
On the other hand, the dimension of any spherical $\th$-twisted conjugacy class in $\tG$ can not be larger than 
$\dim \wt{B}=n \dim B$. Since 
\begin{align*}
(n-1) \dim G - n \dim B &= (n-1) (2\dim B -\dim T) - n \dim B\\
& = (n-3) \dim (B/T) + (\dim (B/T) - \dim T),
\end{align*}
and since $\dim (B/T)$ is the number of positive roots and $\dim T$ is the number of simple roots, 
$\dim \wt{C} > n \dim B$ for every conjugacy class $C$ in $G$ if either $n > 3$ or $n = 3$ and the rank of $G$ is more than $1$, and in
such cases $\wt{C}$ can not be spherical for any $C$.

Let $n = 3$ and $G = SL(2, {\mathbb C})$ or $G=PSL(2, {\mathbb C})$. Then $\dim \wt{C} > 3 \dim B$ if $C$ is not a trivial conjugacy class and thus
$\wt{C}$ can not be spherical. Assume now that $C$ is a trivial conjugacy class consisting of a single central element.
Then $\dim \wt{C} = 2\dim G = 6$. Since $m_{\widetilde{\scriptscriptstyle C}} = (w_0, w_0, w_0)$ and since 
$l(m_{\widetilde{\scriptscriptstyle C}}) + \rk(1-m_{\widetilde{\scriptscriptstyle C}}\tth) = 3+3=6$, it follows from 
from \thref{th-dim-formula} that $\wt{C}$ is spherical.

It remains to consider the case of $n = 2$ and an arbitrary (semi-simple) $G$. In this case, for an arbitrary conjugacy class
$C$ in $G$, $\dim \wt{C} = \dim C + \dim G$ and 
\[
l(m_{\widetilde{\scriptscriptstyle C}}) + \rk(1-m_{\widetilde{\scriptscriptstyle C}}\tth) = 2l(w_0) +\rk(1-(w_0, w_0) \tth).
\]
The linear operator $1-(w_0, w_0) \tth$ on $\h \oplus \h$ is given by $(x, y) \to (x-w_0(y), y-w_0(x))$ and its kernel
is given by $\{(w_0(y), y): y \in \h\}$. Thus $l(m_{\widetilde{\scriptscriptstyle C}}) + \rk(1-m_{\widetilde{\scriptscriptstyle C}}\tth) 
= 2l(w_0)+\dim \h = \dim G$. By \thref{th-dim-formula}, $\wt{C}$ is spherical if and only if $\dim C = 0$, i.e., if and only if
$C$ consists of a single central element.
\end{proof}

\bre{re-referee} Let $n \geq 2$ and let $C$ be a conjugacy class in $G$. For $g \in C$, the 
stabilizer of $(g, e, \ldots, e) \in \wt{C} = \mu_n^{-1}(C)$ in $G^n$ for the $\tilde{\theta}$-twisted conjugation is 
$(Z_g)_{\rm diag} = \{(g, g, \ldots, g): g \in Z_g\}$, where  
$Z_g = \{k \in G: kg = gk\}$ is the centralizer of $g$ in $G$. Thus $\wt{C} \cong G^n/(Z_g)_{\rm diag}$
for any $g \in C$.
When $n = 2$ and $C$ consists of a single central element $g$, $Z_g = G$, and it is well-known that 
$(G \times G)/G_{\rm diag}$ is symmetric and thus spherical. The conclusion in 2) of \leref{le-tth-spherical} can
then be re-interpreted as saying that for a subgroup $Z$ of $G$ that is the centralizer of an element in $G$,
the homogeneous space $(G \times G)/Z_{\rm diag}$ is spherical if and only if $Z = G$. We thank the referee for
pointing this out. 
\ere

\subsection{The Poisson structure $\ptth$ on $G \times G$}
We look at the case of $n = 2$ in more detail. To this end, recall from $\S$\ref{subsec-piG} that the Drinfeld double
of the Poisson Lie group $(G, \pist)$ is the Poisson Lie group  
$(G \times G, \Pist)$, where $G \times G$ has the product Lie group structure and $\Pist$ is the multiplicative Poisson structure 
on $G \times G$ defined by $\Pist = R^r - R^l$ in \eqref{eq-pi-pm}. Recall also from $\S$\ref{subsec-piG} that 
$\Pist$ is invariant under the $T \times T$ given by
\begin{equation}\lb{eq-TT-GG}
(h_1, h_2)\cdot (g_1,g_2) = (h_1g_1h_2^{-1},\, h_1 g_2h_2^{-1}), \hs h_1, h_2 \in T, \, g_1, g_2 \in G. 
\end{equation}

Fix any representative $\dw_0$ of the longest elements $w_0$ of $W$ in $N_T(G)$, and let 
$r_{(\dw_0, \dw_0)}: G \times G \to G \times G$ be the right translation on $G \times G$ by $(\dw_0, \dw_0)$.
Define 
\begin{equation}\lb{eq-tau}
\tau: \;\; G \times G \longrightarrow G \times G, \;\; \tau(g_1, g_2) = (g_1, g_2^{-1}), \hs g_1, g_2 \in G. 
\end{equation}

\bpr{pr-ptth-pim} As Poisson structures on $G \times G$, one has
$\tau(\ptth) = r_{(\dw_0, \dw_0)} (\Pist)$;
\epr

\begin{proof}
We compare the formulas for $\tau(\ptth)$ and $r_{(\dw_0, \dw_0)} (\Pist)$.
By the definition of $\Pist$,
\begin{align*}
(r_{(\dw_0, \dw_0)} (\Pist))(g_1, g_2) &=  r_{(\dw_0, \dw_0)} (\Pist(g_1\dw_0^{-1}, g_2 \dw_0^{-1})) \\
& = l_{(g_1,g_2)} \left(R - \Ad_{(\dw_0^{-1}, \dw_0^{-1})} R\right) + \Pist(g_1, g_2), \hs g_1, g_2 \in G.
\end{align*}
Thus $r_{(\dw_0, \dw_0)} (\Pist)= R^r - \left(\Ad_{(\dw_0^{-1}, \dw_0^{-1})} R\right)^l$.
Let $\{y_i, 1 \leq i \leq k=\dim \h\}$ and $\{E_\al, E_{-\al}: \al \in \Delta_+\}$ be as in $\S$\ref{subsec-dfn}. 
Let $R_\h  = \sum_{i=1}^k (y_i, 0) \wedge (0, y_i) \in \wedge^2(\gog)$ and 
\begin{align*}
R_1 &= \sum_{\al \in \Delta_+} (E_\al, 0) \wedge (0, E_{-\al}),\hs 
R_2 = \sum_{\al \in \Delta_+} (0, E_{\al}) \wedge (E_{-\al}, 0),\\
R_3  &= \frac{1}{2} \sum_{\al \in \Delta_+} \left((E_{\al}, 0) \wedge (E_{-\al}, 0)+ (0, E_\alpha) \wedge (0, E_{-\alpha})\right).
\end{align*}
Using the facts that $\Ad_{\dw_0}$ preserves $\h$ and $\lara_\g$ and that $-\Ad_{\dw_0}$ permutes the positive root spaces, one gets
$R = R_\h + R_1 + R_3$ and $\Ad_{(\dw_0^{-1}, \dw_0^{-1})} R = R_\h - R_2 -R_3.$ Thus
\begin{equation}\lb{eq-dw-dw-pim}
r_{(\dw_0, \dw_0)} (\Pist)= R_\h^r - R_\h^l + (R_1+R_3)^r + (R_2 + R_3)^l.
\end{equation}
On the other hand, applying the explicit formula \eqref{eq-pth-explicit} to $\ptth$, one has
\begin{align*}
\ptth & = \sum_{i=1}^k \tth(y_i, 0)^l \wedge (y_i, 0)^r + \tth(0, y_i)^l \wedge (0, y_i)^r\\
&\;\;\;-\sum_{\al \in \Delta_+} \left((E_\al, 0)^r \wedge \tth(E_{-\al}, 0)^l + (0, E_\al)^r \wedge \tth(0, E_{-\alpha})^l\right)
+R_3^r + R_3^l.
\end{align*}
Using the definitions of $\tth$ and $\tau$ and the fact that $\tau (R_3^r+R_3^l) = R_3^r + R_3^l$, one has
\begin{align*}
\tau(\ptth) & = \sum_{i=1}^k ((y_i, 0)^r \wedge (0, y_i)^r-(y_i, 0)^l \wedge (0, y_i)^l)\\
&\;\;\;+\sum_{\al \in \Delta_+} \left((E_\al, 0)^r \wedge (0, E_{-\al})^r + (0, E_\al)^l \wedge (E_{-\alpha}, 0)^l\right)
+R_3^r + R_3^l,\\
& = R_\h^r -R_\h^l+R_1^r + R_2^l + R_3^r + R_3^l.
\end{align*}
Comparing with \eqref{eq-dw-dw-pim}, one sees that $\tau(\ptth) = r_{(\dw_0, \dw_0)} (\Pist).$
\end{proof}

By \prref{pr-ptth-pim}, one has the $(T \times T)$-equivariant Poisson isomorphism
\[
\phi :=r_{(\dw_0^{-1}, \dw_0^{-1})} \tau: \;\; (G \times G, \ptth)  \lrw (G \times G, \Pist), \;\; 
(g_1, g_2) \longmapsto (g_1\dw_0^{-1}, \, g_2^{-1} \dw_0^{-1}).
\] 
where $T \times T$ acts on $(G \times G, \ptth)$ and on $(G \times G, \Pist)$ respectively by
\begin{align*}
(h_1, h_2)\cdot_{\tth}(g_1, g_2) &= (h_1g_1h_2^{-1}, \, h_2g_2h_1^{-1}),\\
(h_1, h_2)\bullet (g_1, g_2) &= (h_1, w_0(h_2))\cdot (g_1, g_2) = (h_1g_1 w_0(h_2^{-1}), \,h_1g_1 w_0(h_2^{-1})).
\end{align*}
For $u, v \in W$ and a conjugacy class $C$ in $G$, let
\[
G^{u, v}_\sC = \{(k_1, k_2) \in BuB \times B_- v B_-: \, k_1k_2^{-1} \in C\} \subset G \times G.
\]
We will refer to $G^{u, v}_\sC$ the {\it double Bruhat cell associated to $u, v$ and the conjugacy class $C$}.
Note that when $C = \{e\}$, $G^{u, v}_{\sC}$ is isomorphic to the double Bruhat cell 
\[
G^{u, v}= BuB \cap B_-v B_-
\]
which
have been studied intensively \cite{b-f-z:III, f-z:double} in connection  with total positivity and cluster algebras.   
Note that  
\begin{equation}\lb{eq-phi-GuvC}
\phi(\wt{C} \cap (\wt{B} \tilde{w} \wt{B}_-)) = G^{w_1w_0, w_2^{-1}w_0}_{\sC}, \hs w_1, w_2 \in W.
\end{equation} 

\bco{co-pim}
Let $T \times T$ act on $G \times G$ by \eqref{eq-TT-GG}.
The orbits of symplectic leaves of $\Pist$ in $G \times G$ under both $T \times \{e\}$ and
$\{e\} \times T$ 
are the subsets $\GuvC$, where $u, v \in W$ and $C$ a conjugacy class in $G$.
Every $\GuvC$ is a non-empty connected smooth submanifold of $G$ of dimension equal to
$\dim C +  l(u) + l(v) +\dim T$, and the symplectic leaves of $\Pist$ in $\GuvC$ have dimension equal to
$\dim C + l(u) + l(v) + \rk (1-uv^{-1}).$
\eco

\begin{proof}
In the notation of $\S$\ref{subsec-leaves-CBB}, 
$\wt{T}_{\tilde{w}\tth} = \{(h_1w_1(h_2), h_2w_2(h_1)): h_1, h_2 \in T\}$, 
and 
\[
(T \times \{e\})\wt{T}_{\tilde{w}\tth} = 
(\{e\} \times T)\wt{T}_{\tilde{w}\tth} =\wt{T}
\]
 for any $\tilde{w} = (w_1, w_2) \in W \times W$. 
It follows from \thref{th-TT-orb} that under the  action $\cdot_{\tth}$, 
the three tori
$T \times \{e\}$, $\{e\} \times T$, and $T \times T$ have the same orbits of
symplectic leaves of $\ptth$, namely, the subsets
\[
\wt{C} \cap (\wt{B} \tilde{w} \wt{B}_-) = \{(g_1, g_2) \in Bw_1B_- \times B w_2B_-: g_1g_2 \in C\},
\]
where $\tilde{w} = (w_1, w_2) \in W \times W$ and $C$ is a conjugacy class in $G$.
Setting $u = w_1w_0, v = w_2^{-1}w_0$ and using \eqref{eq-phi-GuvC} and \prref{pr-T-orb-ptth}, \coref{co-pim} is proved.
\end{proof}

\bre{re-proof-pim}
\coref{co-pim} is proved by relating $\Pist$ with the Poisson structure $\ptth$ and applying \thref{th-T-orb-mc} to $\ptth$.
A direct proof of \coref{co-pim} is outlined as follows: first, as the Drinfeld double of the Poisson Lie group $(G, \pist)$, 
the symplectic leaves of $\Pist$ in $G \times G$ are the connected components of intersections of 
$(\Gdia, \Gdia)$-double cosets and $(G^*, G^*)$-double cosets (this statement can be proved in a way similar to the proof of
a  statement on the symplectic leaves of the Poisson structure $\pi_+ = R^r + R^l$ on $G \times G$ that can be found in 
\cite{anton-malkin} and \cite[Lemma 6.4]{E-L:grothen}). 
For $u \in W$, define 
\begin{align}\lb{eq-hu}
h_u: &\; \; BuB \longrightarrow T: \; \; \; h_u(n \bar{u} h n^\prime ) = h, \hs n, n^\prime \in N, \, h \in T,\\
\lb{eq-hv}
h_u^\prime: &\; \; B_- u B_- \longrightarrow T: \; \; \; h_u^\prime(n_- \bar{u} h^\prime
 n_-^\prime) = h^\prime,  \hs  n_-, n_-^\prime \in N_-, \, h^\prime \in T.
\end{align}
For $u, v \in W$, let
$T_{u, v} = \{(h^u)^{-1} h^v: \; h \in T\}$, 
where for $h \in T$, $h^u = u^{-1}(h)$. Note that
$\dim T_{u, v} = \rk(1-uv^{-1}).$
It is straightforward to prove that every $(\Gst, \Gst)$-double coset in $G \times G$ is of the form 
$\Gst(\bar{u}h, \, \bar{v})\Gst$ for unique $u, v \in W$ and $h \in T$, and 
\[
\Gst(\bar{u}h, \, \bar{v})\Gst = \{(g_1, g_2) \in BuB \times B_- v B_-: \; \; h_u(g_1) h_v^\prime(g_2) 
\in h T_{u, v}\}.
\]
Using arguments similar to those in the proof of \prref{pr-one-T-leaf}, one can show that 
each $\GuvC = (\Gdia(C, e)\Gdia) \cap (BuB \times B_- v B_-)$ is a smooth connected submanifold of $G \times G$ with the
given dimension and is a $T \times \{e\}$ and a $(\{e\} \times T)$-leaf for $\Pist$ in $G \times G$.
\ere


\begin{thebibliography}{99}

\bibitem{anton-malkin}
A. Alekseev and A. Malkin, {\em Symplectic structures associated to Lie-Poisson groups}, 
Comm. Math. Phys. {\bf 162}(1) (1994), 147 - 173. 

\bibitem{b-f-z:III} A. Berenstein, S. Fomin, and A. Zelevinsky, Cluster Algebras III: upper bounds and Bruhat  cells, {\it Duke Math. J.}
{\bf 126}(1) (2005), 1 - 52.

\bibitem{B:flag} M. Brion, {\it Lectures on the geometry of flag varieties}, {\em Topics in cohomological studies of algebraic varieties},
33-85, Trends in Mathematics, Birkhauser, 2005. 

\bibitem{Can} N. Cantarini, Spherical orbits and quantized enveloping algebras, {\it Comm. Algebra} {\bf 27} (1999), 3439 - 3458.

\bibitem{CCC} N. Cantarini, G. Carnovale, and M. Costantini, Spherical orbits and
representations of ${\mathcal U}_{\epsilon}(\g)$, {\it Trans. Groups} {\bf 10} (1) (2005), 29 - 62.

\bibitem{Car-2008} G. Carnovale, Spherical conjugacy classes and involutions in the Weyl group, {\it Math. Z.} {\bf 260} (1) (2008), 1 - 23.

\bibitem{Car-2012} G. Carnovale, On spherical twisted conjugacy classes, {\it Trans. Groups} {\bf 17} (3) (2012), 615 -5 637.


\bibitem{Chan} K. Y. Chan, {\it Weyl group elements associated to conjugacy classes},
MPhil thesis in Mathematics, The University of Hong Kong, 2010.
 
\bibitem{CLT}
K. Y. Chan, J.-H. Lu, and S. To, On intersections of conjugacy classes and Bruhat cells, {\it Trans. Groups} {\bf 15} (2) (2010), 243 - 260.

\bibitem{C-P:guide}
V. Chari and A. Pressley, A guide to quantum groups, Cambridge University Press, 1994.

\bibitem{Co} M. Costantini, {\it On the coordinate ring of spherical conjugacy classes}, 
Math. Z. {\bf 264} (2010), 327 - 359.
 
\bibitem{Co:bad}
M. Costantini, A classification of unipotent conjugacy classes in bad characteristics, {\it Trans. Amer. Math. Soc.} {\bf 364} (4) (2012), 1997 - 2019.

\bibitem{De} V. Deodhar, On some geometric aspects of Bruhat orderings, I. A finer decomposition of Bruhat cells,
{\it Invent. Math.} {\bf 79} (1985), 499 - 511.

\bibitem{dr:homog}
V. G. Drinfel'd, On Poisson homogeneous spaces of
Poisson-Lie groups, {\em Theo. Math. Phys.}
{\bf 95} (2) (1993), 226 - 227.

\bibitem{En-Et} B. Enriquez and P. Etingof, Quantization of classical dynamical $r$-matrices with nonabelian base, 
{\it Comm. Math. Phys}. {\bf 254} (2005), 603 - 650.

\bibitem{Eting-Schiff}
P. Etingof and O. Schiffmann, {\em Lectures on quantum groups}, 2nd edition, international press, 2002. 

\bibitem{E-L:cplx}
S. Evens and J.-H. Lu, On the variety of Lagrangian
subalgebras, II, {\it Ann. Sci. \'Ecole Norm. Sup}. {\bf{39}} (2) (2006), 347 - 379. 

\bibitem{E-L:grothen}
S. Evens and  J.-H. Lu, Poisson geometry of the Grothendieck resolution of a 
complex semisimple group, {\it Moscow Math. J.} {\bf 7} (4) (special volume 
in honor of V. Ginzburg's 50'th birthday) (2007),  613 - 642.

\bibitem{f-z:double}
S. Fomin and A. Zelevinsky, Double Bruhat Cells and total positivity, {\it J. Amer. Math. Soc.} {\bf{12}} (1999), 335--380.

\bibitem{k-s:quantum}
L. Korogodski and Y. Soibelman, {\em Algebras of functions on
quantum groups, part I}, AMS, Mathematical surveys and
monographs, Vol. 56, 1998.

\bibitem{LiB-M}
D. Li-Bland and E. Meinrenken, Courant algebroids and Poisson geometry, {\it Intern. Math. R. Notices}, 11 (2009) 2106 - 2145.
 
\bibitem{Lu:formula}
J.-H. Lu, On a dimension formula for twisted spherical conjugacy classes in semisimple algebraic groups,
{\it Math. Z.} {\bf 269} (3-4) (2011), 1181 - 1188. 

\bibitem{L-Y:DQ}
J.-H. Lu  and M. Yakimov, Group orbits and regular partitions of Poisson manifolds, {\it Comm. Math. Phys}.
{\bf 283} (3) (2008), 729 - 748.

\bibitem{Panyushev}
D. Panyushev, Complexity and nilpotent orbits. {\it Manuscripta Math.} {\bf  83} (3-4) (1994), 223 - 237.

\bibitem{Ri} R. Richardson, Intersections of double 
cosets in algebraic groups, {\it Indagationes Mathematicae}, Volume 3, Issue 1, (1992),
69 - 77.

\bibitem{Simon-thesis} S. To, On some aspects of a Poisson structure on a complex semisimple Lie group,
PhD thesis in Mathematics, the University of Hong Kong, 2011.

\end{thebibliography}
\end{document}